\documentclass[11pt]{article}
\usepackage{amssymb,amsthm,amsmath,amsthm}
\usepackage[numbers,sort&compress]{natbib}
\newtheorem{theorem.}{\textbf{Theorem}}[section]
\newtheorem{lemma.}{\textbf{Lemma}}[section]
\newtheorem{remark.}{\textbf{Remark}}
\date{}
\begin{document}

\title{\bf  Multiple solutions for superlinear
Klein-Gordon-Maxwell equations\thanks{D.-L. Wu is partially supported
by  NSF of China (No.11801472) and the Youth Science and Technology
Innovation Team of Southwest Petroleum University for Nonlinear
Systems (No.2017CXTD02) and The Science and Technology Innovation Team
of Education Department of Sichuan for Dynamical System and its
Applications (No.18TD0013). H.X. Lin is partially supported by the NSF
of China (No.11701049), the China Postdoctoral Science Foundation
(No.2017M622989), and the Scientific fund of the Education Department
of Sichuan Province (No.18ZB0069).}}
\author{Dong-Lun Wu$^{a,b}$, Hongxia Lin$^{c}$\\
{\small $^{a}$College of Science, Southwest Petroleum University,}\\
{\small Chengdu, Sichuan 610500, P.R. China}\\
{\small $^{b}$Institute of Nonlinear Dynamics, Southwest Petroleum University,}\\
{\small Chengdu, Sichuan 610500, P.R. China}\\
{\small $^{c}$Geomathematics Key Laboratory of Sichuan Province, Chengdu University of
Technology,}\\
{\small Chengdu, Sichuan 610059, P.R. China} }
\date{}
\maketitle
\let\thefootnote\relax\footnotetext{E-mail: wudl2008@163.com; linhongxia18@126.com.}

{\bf Abstract} In this paper, we consider the following
Klein-Gordon-Maxwell equations
\begin{eqnarray*}
\left\{
\begin{array}{ll}
-\Delta u+ V(x)u-(2\omega+\phi)\phi u=f(x,u)+h(x)&\mbox{in $\mathbb{R}^{3}$},\\
-\Delta \phi+ \phi u^2=-\omega u^2&\mbox{in $\mathbb{R}^{3}$},
\end{array}
\right.
\end{eqnarray*}
where $\omega>0$ is a constant, $u$, $\phi : \mathbb{R}^{3}\rightarrow
\mathbb{R}$, $V : \mathbb{R}^{3} \rightarrow\mathbb{R}$ is a potential
function. By assuming the coercive condition on $V$ and some new
superlinear conditions on $f$, we obtain two nontrivial solutions when
$h$ is nonzero and infinitely many solutions when $f$ is odd in $u$
and $h\equiv0$ for above equations.


\section{Introduction and main results}

In this paper, we considered the following Klein-Gordon-Maxwell
equations
\begin{eqnarray}
\left\{
\begin{array}{ll}
-\Delta u+ V(x)u-(2\omega+\phi)\phi u=f(x,u)+h(x)&\mbox{in $\mathbb{R}^{3}$},\\
-\Delta \phi+ \phi u^2=-\omega u^2&\mbox{in $\mathbb{R}^{3}$},
\end{array}
\right.\label{1}
\end{eqnarray}
where $\omega>0$ is a constant, $u$, $\phi$, $V :
\mathbb{R}^{3}\rightarrow \mathbb{R}$. This type of equation has very
interesting physical background which is a model to describe the
nonlinear Klein-Gordon field interacting with the electromagnetic
field. Along with the development of variational methods, many
mathematicians used these methods to investigate the existence and
multiplicity of solutions for Klein-Gordon-Maxwell
equations(see\cite{AP2,1,2,Ca,Carr,Carr2,6,CT,CS,7,19,DM0,8,PL,9,10,11,TZ,12,15,16,22}).
In 2001, V. Benci and D. Fortunato \cite{2} considered the following
systems
\begin{eqnarray}
\left\{
\begin{array}{ll}
-\Delta u+[m^2-(\omega+\phi)^2]u=|u|^{q-2}u&\mbox{in $\mathbb{R}^{3}$},\\
-\Delta \phi+ \phi u^2=-\omega u^2&\mbox{in $\mathbb{R}^{3}$}.
\end{array}
\right.\label{4}
\end{eqnarray}
By using the variational methods, they obtained infinitely many
solitary wave solutions when $|m|>|\omega|$ and $4<q<6$. After this
first work, many mathematicians have treated this problem with variant
cases with $m$, $\omega$ and $q$. In 2004, D'Aprile and Mugnai
\cite{19} dealt with the case $q\in(2,4]$ with
$\displaystyle\sqrt{\left(\frac{q}{2}-1\right)}m>\omega>0$. Some
existence and nonexistence results are obtained for problem (\ref{1})
with variant conditions on $m$, $\omega$ and $q$(see
\cite{19,1,AP2,12,DM0} for more details).

Although we lose the compactness for problem (\ref{1}) since the
problem lies in a unbounded domain, we can also consider this problem
in a radial symmetric space $ H_{r}^{1}\left(\mathbb{R}^{3}\right):=\left\{u \in H^{1}\left(\mathbb{R}^{3}\right)\right.$\\
$\left. : u=u(r),
r=|x|\right\}$. It is known that $
H_{r}^{1}\left(\mathbb{R}^{3}\right)$ is compactly embedded into
$L^{s}\left(\mathbb{R}^{3}\right)$ for $2<s<2^{*}$. If $V(x)$ is not
radial symmetric, there are still some other ways to retrieve the
compactness. One classical way is to assume $V(x)$ to be coercive. In
this paper, we mainly consider the coercive case. This case has also
been studied in many papers(see \cite{3,6,7,8,9,10,11}). Many
authors(see \cite{7,10,13,14}) considered the following coercive
condition.

$(V1)$ Suppose $V \in C(\mathbb{R}^{3},\mathbb{R})$,
$\inf_{\mathbb{R}^{3}}V(x)=V_{0}>0$ and there is a constant $r>0$ such
that
\begin{eqnarray*}
\lim_{|y|\rightarrow+\infty}meas\left(\left\{x\in
 \mathbb{R}^{3}:|x-y|\leq r,\ V(x)\leq M\right\}\right)=0,\ \ \ \forall\ M>0,
\end{eqnarray*}
where $meas(\cdot)$ denote the Lebesgue measure. The condition $(V1)$
was introduced by T. Bartsch et al. [\cite{3}, Lemma 3.1.] to
guarantee the compactness of the embedding. We will use this condition
to prove our theorem.

As we know, the growth of the nonlinear terms is important in showing
the geometric structure of the corresponding functionals(which is
defined in (\ref{3010})) and the boundedness of the Plais-Smale $(PS)$
sequence. In \cite{CT}, the author assumed the following condition.

$(AR)$ there exists a constant $ \theta > 4$ such that
\begin{eqnarray*}
f(x,t)t\geq\theta F(x,t)\ \ \ \mbox{for}\ \ (x,t)\in\ \mathbb{R}^{3}\times\mathbb{R},
\end{eqnarray*}
where $F(x,t)=\int_{0}^{t}f(x,s)ds$. This condition is used to show
the geometric structure of the corresponding functionals and the
boundedness of the $(PS)$ sequences. Hence $(AR)$ has been widely used
in many papers to obtain the existence and multiplicity of elliptic
problem with variational methods. This condition implies that the
growth exponent of the nonlinear terms is more than 4 at infinity. In
order to deal with other nonlinearities with 4-superlinear growth,
many authors(see\cite{8,9}) considered the following growth
conditions.

$(SL1)$ $f(x, t)/t^{3} \rightarrow\infty\ \ \mbox{as}\
|t|\rightarrow+\infty\ \mbox{uniformly in}\ x;$

$(SL2)$ $\frac{1}{4}f(x, t)t-F(x,t) \geq-Dt^{2}\ \ \mbox{for}\ D>0$
and $|t|$ large enough$\ \mbox{uniformly in}\ x.$

$(SL1)$ and $(SL2)$ has been directly or indirectly used to show the
existence and multiplicity of solutions for problem (\ref{1}).
However, above conditions on $f$ eliminate many nonlinearities such as
$F(x,t)=|t|^{5/2}$.

In order to deal with the nonlinearities with growth exponent between
2 and 4, Chen and Song \cite{6} obtained two solutions for problem
(\ref{1}) under $(AR)$ condition just requiring $\theta>2$ with
coercive potentials when $h(x)\not\equiv0$. In 2018, Chen and Tang
\cite{7} introduced the following superlinear condition which is
weaker than $(AR)$ and obtained infinitely many solutions for problem
(\ref{1}).

$(WAR)$ there exist constants $ \theta > 2$ and $K>0$ such that
\begin{eqnarray*}
 f(x,t)t-\theta F(x,t)+Kt^{2}\geq0\ \ \ \mbox{for}\ \ (x,t)\in\ \mathbb{R}^{3}\times\mathbb{R}.
\end{eqnarray*}

There are still many functions cannot be included in above conditions.
In this paper, we introduce some new superlinear conditions and an
example is given to show the difference from previous conditions. Now,
we state our main results.

\begin{theorem.}\label{th1.1} Suppose that $V$ satisfies $(V1)$ and $f \in C(\mathbb{R}^{3} \times \mathbb{R},\mathbb{R})$
satisfies

$(F1)$ $f(x,t)=o(|t|)$ as $t\rightarrow0$ for all $x\in
\mathbb{R}^{3}$;

$(F2)$ $F(x,0)=0$ for all $x\in\mathbb{R}^{3}$ and there exist
$\tau\in\left(2,6\right)$ and $d_{1}>0$ such that
\begin{eqnarray*}
|f(x,t)|\leq d_{1}(|t|+|t|^{\tau-1})\ \ \ \ \mbox{for all}\ \
(x,t)\in \mathbb{R}^{3}\times\mathbb{R};
\end{eqnarray*}

$(F3)$ let $\widetilde{F}(x,t)=f(x,t)t-2 F(x,t)$, then there exist
$d_{2}>0$ and $r_{0}>0$ such that
\begin{eqnarray*}
\widetilde{F}(x,t)\geq d_{2}|t|^{\tau}\ \  \mbox{for all}\ \ |t|\geq r_{0} \ \ \mbox{and}\  \ x\in
\mathbb{R}^{3};
\end{eqnarray*}

$(F4)$ $F(x,t) / t^{2} \rightarrow +\infty$ as $|t| \rightarrow
\infty$ uniformly in $x$;

$(F5)$ $f(x,t)t\geq 2F(x,t)\geq0$ for all
$(x,t)\in\mathbb{R}^{3}\times\mathbb{R}$.

Then there is a constant $m>0$ such that, for any $h\not\equiv0$
satisfying $\|h\|_{2}<m$, problem (\ref{1}) possesses at least two
nontrivial solutions.
\end{theorem.}

\begin{theorem.}\label{th1.2} Suppose that $(V1)$,
$(F2)$-$(F5)$ hold, $h(x)\equiv0$ and $f(x,-t)=-f(x,t)$, for all
$(x,t)\in\mathbb{R}^{3}\times\mathbb{R}$. Then problem (\ref{1})
possesses infinitely many solutions.
\end{theorem.}

\begin{remark.}
In 2018, Shi and Chen\cite{11} obtained two solutions for problem
(\ref{1}) with nonzero perturbation by using the combination of a
cut-off functional and a Pohozaev type identity. Although $f(x,t)$ was
required to satisfy some very weak growth conditions in their paper,
they needed some smooth conditions on the gradients of $V(x)$ and
$h(x)$, which are not needed in our theorems.
\end{remark.}

\begin{remark.}Setting $4>p>2$, $0<\epsilon<p-2$, consider
\begin{eqnarray}
F(x,t)=|t|^{p}+a(p-2)|t|^{p-\epsilon}\sin^{2}(|t|^{\epsilon}/\epsilon).\label{50}
\end{eqnarray}
For any $\theta>2$ and $K>0$, let
$\max\left\{0,\frac{p-\theta}{p-2}\right\}<a<1$ and
$t_{n}=\left(\epsilon\left(n\pi+\frac{3\pi}{4}\right)\right)^{1/\epsilon}$,
then
\begin{eqnarray*}
&&f(x,t_{n})t_{n}-\theta F(x,t_{n})-Kt_{n}^{2}\\
&=&(p-\theta)|t_{n}|^{p}+a(p-2)(p-\theta-\epsilon)|t_{n}|^{p-\epsilon}\sin^{2}(|t_{n}|^{\epsilon}/\epsilon)+a(p-2)|t_{n}|^{p}\sin(2|t_{n}|^{\epsilon}/\epsilon)-Kt_{n}^{2}\\
&=&|t_{n}|^{p}\left[(p-\theta)-a(p-2)+\frac{a(p-2)(p-\theta-\epsilon)}{2|t_{n}|^{\epsilon}}\right]-Kt_{n}^{2}\\
&\rightarrow&-\infty\ \ \mbox{as}\ \ n\rightarrow\infty.
\end{eqnarray*}
Hence (\ref{50}) does not satisfy $(WAR)$ or the following condition

$(FSL)$ $F(x, t)/t^{4} \rightarrow+\infty\ \ \mbox{as}\
|t|\rightarrow+\infty\ \mbox{uniformly in}\ x.$

However, it is easy to see that, for all $x\in \mathbb{R}^{3}$,
\begin{eqnarray*}
f(x,t)t-2 F(x,t)
&=&(p-2)|t|^{p}\left[\left(1+a\sin(2|t|^{\epsilon}/\epsilon)\right)+\frac{a(p-2-\epsilon)\sin^{2}(|t|^{\epsilon}/\epsilon)}{|t|^{\epsilon}}\right]\\
&\geq&\frac{1}{2}(1-a)(p-2)|t|^{p}
\end{eqnarray*}
for $|t|$ large enough. Then we can check that (\ref{50}) satisfies
conditions $(F1)$-$(F5)$.
\end{remark.}

\section{Preliminaries}

Let
$$E:=\left\{u\in
H^1({\mathbb{R}}^3):\int_{{\mathbb{R}}^3}\left(|\nabla
u|^2+V(x)u^2\right)dx<\infty \right\}.$$
Then $E$ is a Hilbert space
with the inner product
$$(u,v)_E=\int_{{\mathbb{R}}^3}\left(\nabla
u\cdot \nabla v+V(x)uv\right)dx$$
and the norm $\|u\|_E=(u,u)_E^{1/
2}$. Define the function space
$$D^{1,2}({\mathbb{R}}^3):=\left\{u\in
L^{6}({\mathbb{R}}^3): |\nabla u| \in L^2({\mathbb{R}}^3) \right\}$$
with the norm
$$\|u\|_{D^{1,2}}:=\left(\int_{{\mathbb{R}}^3}|\nabla
u|^2dx\right)^{\frac{1}{2}}.$$ Under $(V1)$, the embedding
$E\hookrightarrow L^s({{\mathbb{R}}^3})$ are compact for any $s\in [2,
6)$(see \cite{3}). Hence for each $s\in [2,6)$, there exists a
constant $C_{s}>0$ such that
\begin{equation}\label{2}
\|u\|_{L^s}\leq C_{s} \|u\|_{E},\ \ \ \ \forall u\in E.
\end{equation}
Obviously, problem (\ref{1}) has a variational structure. Consider
 $J: E \times D^{1,2}({\mathbb{R}}^3)\rightarrow \mathbb{R}$
defined by $$J(u,\phi)=\frac{1}{2} \int_{{\mathbb{R}}^3}\left(|\nabla
u|^2-|\nabla
\phi|^2+[V(x)-(2\omega+\phi)\phi]u^2\right)dx-\int_{{\mathbb{R}}^3}F(x,u)dx-\int_{{\mathbb{R}}^3}h(x)udx.
$$ Evidently, $J$ belongs to $C^1(E\times
D^{1,2}({\mathbb{R}}^3),\mathbb{R})$. Thus,  the pair $(u, \phi)$ is a
weak solution of problem (\ref{1})
 if and only if it is a critical point of $J$ in $E\times
 D^{1,2}({\mathbb{R}}^3)$. We can also see that $J$ is strong
 indefinite. To reduce this functional, we need the following lemmas.

\begin{lemma.} (see\cite{19})\label{le2.2} For
every $u\in H^1({\mathbb{R}}^3)$ there
 exists a unique $\phi _u\in D^{1,2}({\mathbb{R}}^3)$
 which solves the second equation of problem (\ref{1}).  Furthermore,

(i) $-\omega \leq \phi_u \leq 0$ in $\mathbb{R}^3$;

(ii) if $u$ is radially symmetric,  $\phi_u$ is radial too.
\end{lemma.}
We can consider the functional $I:E \rightarrow {\mathbb{R}}$ defined
by $I(u)=J(u,\phi_u)$.  Therefore
\begin{equation}\label{3010}
I(u)=\frac{1}{2} \int_{{\mathbb{R}}^3}\left(|\nabla
u|^2+V(x)u^2-\omega\phi_u
u^2\right)dx-\int_{{\mathbb{R}}^3}F(x,u)dx-\int_{{\mathbb{R}}^3}h(x)udx.
\end{equation}
and we have, for any $u, v \in E$,
\begin{equation}\label{301}
\langle I'(u),v\rangle=\int_{{\mathbb{R}}^3}\left(\nabla u \cdot \nabla v
+V(x)uv-(2\omega+\phi_u)\phi_u
uv-f(x,u)v-h(x)v\right)dx.
\end{equation}

\begin{lemma.}\label{le2.3} The following statements
 are equivalent:

(1) $(u,\phi)\in H^1({\mathbb{R}}^3) \times D^{1,2}({\mathbb{R}}^3)$
is a critical point of $J$ (i.e. $(u,\phi)$ is a solution of system
(\ref{1}) );

(2) $u$ is a  critical point of $I$ and $\phi=\phi_u$.
\end{lemma.}
By Lemmas (\ref{le2.2}) and (\ref{le2.3}), we only need to look for
the critical points of $I$ to show the existence and multiplicity of
critical points for $J$. Some details can be found in \cite{10}.
Subsequently, we prove our theorems with mountain pass theorem and a
abstract critical point theorem introduced by T. Bartsch.

\section{Proof of Theorem \ref{th1.1}}

\begin{lemma.}\label{le2.1}
Suppose that $(V1)$, $(F1)$ and $(F2)$ hold. Let $h\in
 L^2(\mathbb{R}^{3})$, then there exist some constants $\rho, \alpha, m>0$ such
 that $I(u)\big |_{\|u\|_E = \rho}\geq \alpha$ for all $h$ satisfying $\|h\|_{L^2}<
 m$.
\end{lemma.}

{\bf Proof.} The proof is similar to Lemma 3.1 in
\cite{6}.~~~~~~~~~~~~~~~~~~~~~~~~~~~~~~~~~~~~~~~~~~~~~~~~~~~~~~~~~~~~~~~~~~~~~~~~~~~$\Box$

\begin{lemma.}
Suppose that $(V1)$, $(F4)$, $(F5)$ hold, then there exists $v\in E$
with $\|v\|_E>\rho$ such that $I(v)<0$, where $\rho$ is given in Lemma
\ref{le2.1}.
\end{lemma.}

{\bf Proof.} Set $e \in C^{\infty}_{0}(B_{1}(0),\mathbb{R})$, where
$B_{1}(0)=\{x\in \mathbb{R}^{3}:|x|\leq1\}$, such that $\|e\|_{E}=1$
and
$A=\frac{2(1+\int_{B_{1}(0)}\omega^{2}e^{2}dx)}{\int_{B_{1}(0)}e^{2}dx}$.
It follows from $(F4)$ and $(F5)$ that there exists $Q>0$ such that
\begin{eqnarray}
 F(x,t)\geq A(t^{2}-Q^{2})\label{33}
\end{eqnarray}
for all $(x,t)\in B_{1}(0)\times \mathbb{R}$.  By (\ref{33}) and Lemma
\ref{le2.2}, for every $\eta \in \mathbb{R}^{+}$, we have
\begin{eqnarray*}
I(\eta e)&=&\frac{\eta^2}{2}
\int_{B_{1}(0)}\left(|\nabla e|^2+V(x)e^2-\omega\phi_{\eta e}
e^2\right)dx-\int_{B_{1}(0)}F(x,\eta e)dx-\eta\int_{B_{1}(0)}h(x)edx\\
&\leq&
\frac{\eta^{2}}{2}\left(1+\int_{B_{1}(0)}\omega^{2}e^{2}dx-2A\int_{B_{1}(0)}e^{2}dx\right)+\frac{4\pi AQ^{2}}{3}-\eta\int_{B_{1}(0)}h(x)edx\nonumber\\
&\leq& \mbox{}
-\frac{3\eta^{2}}{2}\left(1+\int_{B_{1}(0)}\omega^{2}e^{2}dx\right)+\frac{4\pi AQ^{2}}{3}-\eta\int_{B_{1}(0)}h(x)edx.
\end{eqnarray*}
Then $I(\eta e)\rightarrow-\infty$ as $\eta\rightarrow+\infty$.
Therefore, there exists $\eta_{0}>0$ such that $I(\eta_{0}e)<0$ and
$\|\eta_{0}e\|_{E}> \varrho$. Let $v=\eta_{0}e$, we can see $I(v)<0$,
which proves this
lemma.~~~~~~~~~~~~~~~~~~~~~~~~~~~~~~~~~~~~~~~~~~~~~~~~~~~~~~~~~~~~~~~~~~~~~~~~~$\Box$

\begin{lemma.}
Assume that $(V1)$ and $(F2)-(F5)$ hold,   then $I$ satisfies the
$(PS)$ condition.
\end{lemma.}

{\bf Proof.} Suppose that $\{u_{n}\}\subset E$ is a sequence such that
$\{I(u_{n})\}$ is bounded and $I'(u_{n}) \rightarrow 0$ as $n
\rightarrow \infty$. Then there exists a constant $\overline{M} > 0$
such that
\begin{eqnarray}
|I(u_{n})| \leq \overline{M},\ \ \ \ \
\|I'(u_{n})\| \leq \overline{M}.\label{47}
\end{eqnarray}
Now we prove that $\{u_{n}\}$ is bounded in $E$. Arguing in an
indirect way, we assume that $\|u_{n}\|_{E}\rightarrow+\infty$ as
$n\rightarrow \infty$. Set $w_{n}=\frac{u_{n}}{\|u_{n}\|_{E}}$. Then
$\|w_{n}\|_{E}=1$ and there exists a subsequence of $\{w_{n}\}$, still
denoted by $\{w_{n}\}$, such that $w_{n}\rightharpoonup w_{0}$ in $E$.
Then we have
\begin{eqnarray*}
w_{n}\rightarrow w_{0}\ \ \mbox{in}\ L^{s}(\mathbb{R}^{3})\ \ \ \mbox{for any}\ \ s\in[2,6).\label{20}
\end{eqnarray*}
Let $\Omega=\{x\in \mathbb{R}^{3}|\ |w_{0}(x)|>0 \}$. If
$meas(\Omega)>0$, we have $|u_{n}|\rightarrow+\infty$ as $n\rightarrow
\infty$ for a.e. $x\in \Omega$. On one hand, by $(F4)$, $(F5)$ and
Fatou's Lemma, we obtain
\begin{eqnarray*}
\liminf_{n\rightarrow \infty}\int_{\mathbb{R}^3}\frac{F(x,u_{n})}{\|u_{n}\|_{E}^{2}}dx&\geq&\liminf_{n\rightarrow \infty}\int_{\Omega}\frac{F(x,u_{n})}{\|u_{n}\|_{E}^{2}}dx{\nonumber}\\
&=& \mbox{}\liminf_{n\rightarrow \infty}\int_{\Omega}\frac{F(x,u_{n})}{|u_{n}|^{2}}|w_{n}|^{2}dx{\nonumber}\\
&=& \mbox{}+\infty.
\end{eqnarray*}
On the other hand, by (\ref{47}) and Lemma \ref{le2.2}, we get
\begin{eqnarray*}
\left|\int_{\mathbb{R}^{3}}\frac{F(x,u_{n})}{\|u_{n}\|_{E}^{2}}dx-\frac{1}{2}\right|&=&\left|-\frac{I(u_{n})}{\|u_{n}\|_{E}^{2}}-\frac{\omega}{2}\int_{\mathbb{R}^{3}}\frac{\phi_u
u_{n}^2}{\|u_{n}\|_{E}^{2}}dx-\int_{\mathbb{R}^{3}}\frac{h(x)u_{n}}{\|u_{n}\|_{E}^{2}}dx\right|\nonumber\\
&\leq&\frac{\overline{M}}{\|u_{n}\|_{E}^{2}}+\frac{\omega^{2}C_{2}}{2}+\frac{C_{2}\|h\|_{2}}{\|u_{n}\|_{E}}
\rightarrow\frac{\omega^{2}C_{2}}{2}\ \ \ \mbox{as}\ \ \ n\rightarrow\infty,
\end{eqnarray*}
which is a contradiction. Then we have $meas(\Omega)=0$, which implies
that $w_{0}=0$ a.e. $x\in \mathbb{R}^{3}$ and $w_{n}\rightarrow0$ in
$L^{s}(\mathbb{R}^{3})$($2\leq s<6$). By~(\ref{3010}), (\ref{301}),
$(F3)$ and $(F5)$, we obtain
\begin{eqnarray*}
2I(u_{n})-\langle I'(u_{n}),u_{n}\rangle
&=& \mbox{}\int_{\mathbb{R}^{3}}(\omega+\phi_{u_{n}})\phi_{u_{n}}
u_{n}^2dx+ \int_{\mathbb{R}^{3}}\widetilde{F}(x,u_{n})dx-\int_{\mathbb{R}^{3}}h(x)u_{n}dx{\nonumber}\\
&\geq& \mbox{}\int_{\mathbb{R}^{3}}(\omega+\phi_{u_{n}})\phi_{u_{n}}
u_{n}^2dx+ d_{2}\int_{|u_{n}|\geq r_{0}}|u_{n}|^{\tau}dx-\int_{\mathbb{R}^{3}}h(x)u_{n}dx,
\end{eqnarray*}
which implies that
\begin{eqnarray}\label{3}
\lim_{n\rightarrow\infty}\frac{\int_{|u_{n}|\geq r_{0}}|u_{n}|^{\tau}dx}{\|u_{n}\|_{E}^{2}}=0.
\end{eqnarray}
By (\ref{301}), (\ref{3}), $(F2)$ and (\ref{2}), one sees that
\begin{eqnarray*}
\frac{\overline{M}+1}{\|u_{n}\|_{E}}&\geq&\frac{\langle I'(u_{n}),u_{n}\rangle}{\|u_{n}\|_{E}^{2}}\\
&\geq& 1- d_{1}\int_{\mathbb{R}^{3}}\left( w_{n}^{2}+\frac{|u_{n}|^{\tau}}{\|u_{n}\|_{E}^{2}}\right)dx-\frac{\|h\|_{2}\|u_{n}\|_{2}}{\|u_{n}\|_{E}^{2}}{\nonumber}\\
&\geq& 1- d_{1}\int_{\mathbb{R}^{3}}w_{n}^{2}dx-d_{1}\left(\frac{\int_{|u_{n}|\geq r_{0}}|u_{n}|^{\tau}dx}{\|u_{n}\|_{E}^{2}}+r_{0}^{\tau-2}\int_{|u_{n}|\leq r_{0}}w_{n}^{2}dx\right)-\frac{C_{2}\|h\|_{2}}{\|u_{n}\|_{E}}{\nonumber}\\
&\rightarrow&1\ \ \ \mbox{as}\ \ \ n\rightarrow\infty,
\end{eqnarray*}
which  is a contradiction, then $\{u_{n}\}$ is bounded in $E$. Similar
to Lemma 3.3 in \cite{6}, we can see that $\{u_{n}\}$ has a strong
convergent subsequence. Then $I$ satisfies the $(PS)$
condition.~~~~~~~~~~~~~~~~~~~~~~~~~~~~~~~~~~~~~~~~~~~~~~~$\Box$

{\bf Proof of Theorem 1.1.} In order to obtain two nontrivial
solutions for problem (\ref{1}), we will apply the Ekeland's
variational principle and the mountain pass theorem to the functional
$I$. The rest proof of Theorem 1.1 is similar to the proof of Theorem
1.2 in \cite{6}.

\section{Proof of Theorem \ref{th1.2}}

In order to obtain infinitely many solutions of (\ref{1}), similar to
\cite{10}, we shall use the abstract critical point theorem introduced
by T. Bartsch in \cite{17}. Let space $X$ be reflexive and separable,
then there exist ${e_{i}} \in X$ and ${f_{i}} \in X^{*}$ such that $X
= \overline{\langle e_{i}, i \in \mathbb{N}\rangle}$, $X^{*} =
\overline{\langle f_{i}, i \in \mathbb{N}\rangle}$, $\langle e_{i},
f_{j}\rangle = \delta_{i,j}$, where $\delta_{i,j}$ denotes the
Kronecker symbol. Subsequently, put $$X_{k} = span\{e_{k}\},\ \ \
Y_{k} = \bigoplus_{i=1}^{k}  X_{i},\ \ \  Z_{k} =
\overline{\bigoplus_{i=k}^{\infty} X_{i}}.$$

Now we state the critical points theorem by T. Bartsch.
\begin{lemma.}\label{le2.4} Assume $\phi \in C^{1}(X, \mathbb{R})$ satisfies the $(PS)$ condition, $\phi(-u)=\phi(u)$. For every $k \in \mathbb{N}$, there exists
$\rho_{k}>r_{k}>0$, such that

(i) $a_{k}:=\max_{u\in Y_{k}, \|u\|=\rho_{k}} \phi(u) \leq 0$;

(ii) $b_{k}:=\inf_{u\in Z_{k}, \|u\|=r_{k}} \phi(u)\rightarrow
+\infty$ as $k \rightarrow\infty$.

Then $\phi$ has a sequence of critical values tending to $+\infty$
\end{lemma.}

\begin{lemma.} Assume that $(V1)$, $(F2)$, $(F4)$ and $(F5)$ hold. Then for every $k \in \mathbb{N}$,
there exists $\rho_{k} > r_{k} > 0$, such that (i) $a_{k} :=
\inf_{u\in Z_{k},\|u\|_{E}=r_{k}} I(u) \rightarrow+\infty$ as $k
\rightarrow\infty$; (ii) $b_{k} := \max_{u\in
Y_{k},\|u\|_{E}=\rho_{k}} I(u) \leq 0$.
\end{lemma.}

{\bf Proof.} It follows from $(F2)$ that
\begin{eqnarray}
|F(x,t)|\leq d_{1}\left(\frac{1}{2} t^{2}+\frac{1}{\tau}|t|^{\tau}\right),\ \ \ \ \
 \forall (x,t) \in \mathbb{R}^{3}\times\mathbb{R}.\label{28}
\end{eqnarray}
For any $k\in \mathbb{N}$ and $p\in[2,6)$, we set $$
\beta_{k}(p)=\sup_{u\in Z_{k},\ \|u\|_{E}=1}\|u\|_{p}. $$ Similar with
Lemma 2.8 in \cite{10}, we have $\beta_{k}(p)\rightarrow0$ as
$k\rightarrow\infty$. Letting $
r_{k}=\left(\frac{\tau}{8d_{1}\beta_{k}^{\tau}(\tau)}\right)^{\frac{1}{\tau-2}}$,
for any $u\in Z_{k}$, it follows from (\ref{3010}), Lemma \ref{le2.2}
and (\ref{28}) that
\begin{eqnarray*}
I(u)&=&
\frac{1}{2}
\int_{\mathbb{R}^{3}}\left(|\nabla u|^2+V(x)u^2-\omega\phi_{u}
u^2\right)dx-\int_{\mathbb{R}^{3}}F(x, u)dx{\nonumber}\\
&\geq& \mbox{}\frac{1}{2}\|u\|_{E}^{2}-d_{1}\left(\frac{\beta_{k}^{2}(2)}{2}\|u\|_{E}^{2}+\frac{\beta_{k}^{\tau}(\tau)}{\tau}\|u\|_{E}^{\tau}\right)\\
&\geq& \mbox{}\frac{1}{8}\|u\|_{E}^{2}
\end{eqnarray*}
for $k$ large enough. Thus we obtain $ a_{k}=\inf_{u\in
Z_{k},\|u\|_{E}=r_{k}}I(u)\geq\frac{1}{8}r_{k}^{2}\rightarrow+\infty\
\ \mbox{as}\ \ k\rightarrow\infty$.

Subsequently, for any $u\in Y_{k}$ and $\delta>0$, set
$$\Gamma_{\delta}(u)=\{x\in \mathbb{R}^{3}:\ |u|\geq\delta \|u\|_{E}\}.$$
Similar to Lemma 2.6 in \cite{18}, there exists $\varepsilon_{1}>0$
such that $$meas(\Gamma_{\varepsilon_{1}}(u))\geq\varepsilon_{1}.$$
From $(F4)$, there exists $r_{\infty}>0$ such that
\begin{eqnarray}
F(x,u)\geq \frac{(1+\omega^{2}C_{2}^{2})}{\varepsilon_{1}^{3}}|u|^{2}\geq
\frac{(1+\omega^{2}C_{2}^{2})}{\varepsilon_{1}}\|u\|_{E}^{2}\label{12}
\end{eqnarray}
for all $u\in Y_{k}$ and $x\in \Gamma_{\varepsilon_{1}}(u)$ with
$\|u\|_{E}\geq \frac{r_{\infty}}{\varepsilon_{1}}$. We can choose
$\rho_{k}>\max\left\{\frac{r_{\infty}}{\varepsilon_{1}},r_{k}\right\}$,
then for any $u\in Y_{k}$ with $\|u\|_{E}=\rho_{k}$, it follows from
(\ref{3010}), $(F4)$, Lemma \ref{le2.2}, (\ref{12}) and $(F5)$ that
\begin{eqnarray*}
I(u)&=&
\frac{1}{2}
\int_{\mathbb{R}^{3}}\left(|\nabla u|^2+V(x)u^2-\omega\phi_{u}
u^2\right)dx-\int_{\mathbb{R}^{3}}F(x, u)dx{\nonumber}\\
&\leq& \mbox{}\frac{1}{2}\|u\|_{E}^{2}+\frac{\omega^{2}}{2}\|u\|_{2}^{2}-\int_{\Gamma_{\varepsilon_{1}}(u)}F(x,u)dx{\nonumber}\\
&\leq& \mbox{}-\frac{(1+\omega^{2}C_{2}^{2})}{2}\|u\|_{E}^{2},
\end{eqnarray*}
which means $b_{k}\leq0$ for $\rho_{k}$ large enough. Then we finish
the proof of this lemma.~~~~~~~~~~~~~~~~~~~~~~~~~~~~~~~~~~~~~~~~$\Box$

{\bf Proof of Theorem 1.2.} Similar to the proof of Theorem
\ref{th1.1}, we see that $I$ satisfies the $(PS)$ condition.
Furthermore, $I(-u)=I(u)$, then we obtain our conclusion by using the
Lemma \ref{le2.4}.

\section{Acknowledgments}

The authors are grateful to the referees for the helpful comments
which improve the writing of the paper. This paper was finished when
D.-L. Wu was visiting Utah State University with the support of China
Scholarship Council(No.201708515186); he is grateful to the members in
the Department of Mathematics and Statistics at Utah State University
for their invitation and hospitality.

\def\refname{References}


\begin{thebibliography}{99}

\bibitem[{Azzollini and Pomponio(2010)}]{AP2}
    \bibinfo{author}{A.~Azzollini}, \bibinfo{author}{A.~Pomponio},
  \bibinfo{title}{Ground state solutions for the nonlinear
  {K}lein-{G}ordon-{M}axwell equations}, \bibinfo{journal}{Topol. Methods
  Nonlinear Anal.} \bibinfo{volume}{35}~(\bibinfo{number}{1})
  (\bibinfo{year}{2010}) \bibinfo{pages}{33--42}.

\bibitem{1} A. Azzollini, L. Pisani, A. Pomponio, Improved estimates
    and a limit case for the electrostatic Klein-Gordon-Maxwell system,
    Proc. Roy. Soc. Edinburgh Sect. A 141 (3) (2011) 449-463.


\bibitem{2} V. Benci, D. Fortunato, The nonlinear Klein-Gordon
    equation
    coupled with the Maxwell equations, Nonlinear Anal. 47 (2001)
    6065-6072.

\bibitem{17} T. Bartsch, Infinitely many solutions of a symmetric
    Dirichlet problem, Nonlinear Anal. 20 (10) (1993) 1205-1216.

\bibitem{3} T. Bartsch, Z.-Q. Wang, M. Willem, The Dirichlet problem
    for superlinear elliptic equations, in: Stationary Partial
    Differential Equations. Vol. II, in: Handb. Differ. Equ.,
    Elsevier/North-Holland, Amsterdam, 2005, pp. 1-55.

\bibitem[{Cassani(2004)}]{Ca} \bibinfo{author}{D.~Cassani},
    \bibinfo{title}{Existence and non-existence of
  solitary waves for the critical {K}lein-{G}ordon equation coupled with
  {M}axwell's equations}, \bibinfo{journal}{Nonlinear Anal.}
  \bibinfo{volume}{58}~(\bibinfo{number}{7-8}) (\bibinfo{year}{2004})
  \bibinfo{pages}{733--747}.

\bibitem[{Carri{\~a}o et~al.(2011)Carri{\~a}o, Cunha, and
    Miyagaki}]{Carr} \bibinfo{author}{P.~C. Carri{\~a}o},
    \bibinfo{author}{P.~L. Cunha},
  \bibinfo{author}{O.~H. Miyagaki}, \bibinfo{title}{Existence results for the
  {K}lein-{G}ordon-{M}axwell equations in higher dimensions with critical
  exponents}, \bibinfo{journal}{Commun. Pure Appl. Anal.}
  \bibinfo{volume}{10}~(\bibinfo{number}{2}) (\bibinfo{year}{2011})
  \bibinfo{pages}{709--718}.

\bibitem[{Carri{\~a}o et~al.(2012)Carri{\~a}o, Cunha, and
    Miyagaki}]{Carr2} \bibinfo{author}{P.~C. Carri{\~a}o},
    \bibinfo{author}{P.~L. Cunha},
  \bibinfo{author}{O.~H. Miyagaki}, \bibinfo{title}{Positive ground state
  solutions for the critical {K}lein-{G}ordon-{M}axwell system with
  potentials}, \bibinfo{journal}{Nonlinear Anal.}
  \bibinfo{volume}{75}~(\bibinfo{number}{10}) (\bibinfo{year}{2012})
  \bibinfo{pages}{4068--4078}.

\bibitem{6} S.J. Chen, S.Z. Song, Multiple solutions for
    nonhomogeneous Klein-Gordon-Maxwell equations on $\mathbb{R}^{3}$, Nonlinear
    Analysis: Real World Applications 22 (2015) 259-271.

\bibitem[{Chen, Shang-Jie and Tang, Chun-Lei(2010)}]{CT}
    \bibinfo{author}{S.-J.~Chen}, \bibinfo{author}{C.-L.~Tang},
   \bibinfo{title}{Multiple solutions for nonhomogeneous
              {S}chr\"odinger-{M}axwell and {K}lein-{G}ordon-{M}axwell
              equations on {$\mathbb{R}\sp 3$}}, \bibinfo{journal}{NoDEA Nonlinear Differential Equations Appl.}
  \bibinfo{volume}{17}~(\bibinfo{number}{5}) (\bibinfo{year}{2010})
  \bibinfo{pages}{559-574}.

\bibitem[{Candela and Salvatore(2009)}]{CS} \bibinfo{author}{A.~M.
    Candela}, \bibinfo{author}{A.~Salvatore},
  \bibinfo{title}{Multiple solitary waves for non-homogeneous
  {K}lein-{G}ordon-{M}axwell equations}, in: \bibinfo{booktitle}{More Progresses in Analysis
Proceedings of the 5th International ISAAC Congress ({C}atania, 2005)
},  \bibinfo{pages}{753--762},   \bibinfo{year}{2009}.

\bibitem{7} S.T. Chen, X.H. Tang, Infinitely many solutions and least
    energy solutions for Klein-Gordon-Maxwell systems with general superlinear
    nonlinearity, Comput. Math. Appl. 75 (2018) 3358-3366.


\bibitem{19} T. D'Aprile, D. Mugnai, Solitary waves for nonlinear
    Klein-Gordon-Maxwell and Schr\"{o}dinger-Maxwell equations, Proc. Roy.
    Soc. Edinburgh Sect. A 134 (5) (2004) 893-906

\bibitem[{D'Aprile and Mugnai(2004{\natexlab{b}})}]{DM0}
    \bibinfo{author}{T.~D'Aprile}, \bibinfo{author}{D.~Mugnai},
  \bibinfo{title}{Non-existence results for the coupled
  {K}lein-{G}ordon-{M}axwell equations}, \bibinfo{journal}{Adv. Nonlinear
  Stud.} \bibinfo{volume}{4}~(\bibinfo{number}{3})
  (\bibinfo{year}{2004}{\natexlab{b}}) \bibinfo{pages}{307--322}.

\bibitem{8} L. Ding, L. Li, Infinitely many standing wave solutions
    for
    the nonlinear Klein-Gordon-Maxwell system with sign-changing
    potential, Comput. Math. Appl. 68 (2014)
    589-595.

\bibitem[{d'Avenia and Pisani(2002)}]{PL}
    \bibinfo{author}{P.~d'Avenia}, \bibinfo{author}{L.~Pisani},
  \bibinfo{title}{Nonlinear {K}lein-{G}ordon equations coupled with
  {B}orn-{I}nfeld type equations}, \bibinfo{journal}{Electron. J. Differential
  Equations} \bibinfo{pages}{26} (\bibinfo{year}{2002}) 1-13.

\bibitem{9} X.M. He, Multiplicity of solutions for a nonlinear
    Klein-Gordon-Maxwell system, Acta Appl. Math. 130 (2014) 237-250.


\bibitem{10} L. Li, C.-L. Tang, Infinitely many solutions for a
    nonlinear Klein-Gordon-Maxwell system, Nonlinear Anal. 110 (2014) 157-169.



\bibitem{11} H.X. Shi, H.B. Chen, Multiple positive solutions for
    nonhomogeneous Klein-Gordon-Maxwell equations, Appl. Math. Comput. 337 (2018)
    504-513.

\bibitem[{Teng and Zhang(2011)}]{TZ} \bibinfo{author}{K.~Teng},
    \bibinfo{author}{K.~Zhang},
  \bibinfo{title}{Existence of solitary wave solutions for the nonlinear
  {K}lein-{G}ordon equation coupled with {B}orn-{I}nfeld theory with critical
  {S}obolev exponent}, \bibinfo{journal}{Nonlinear Anal.}
  \bibinfo{volume}{74}~(\bibinfo{number}{12}) (\bibinfo{year}{2011})
  \bibinfo{pages}{4241--4251}.

\bibitem{12} F.Z. Wang; Ground state solutions for the electrostatic
    nonlinear Klein-Gordon-Maxwell system, Nonlinear Anal. 74 (14)
    (2011), 4796-4803.


\bibitem{13} D.-L. Wu, F. Li, Solutions for fourth-order Kirchhoff type elliptic equations involving concave-convex nonlinearities in $\mathbb{R}^{N}$, Comput. Math. Appl. 79(2)(2020), 489-499.

\bibitem{14} D.-L. Wu, X. Yu, New homoclinic orbits for Hamiltonian systems with asymptotically quadratic growth at infinity, Qual. Theory Dyn. Syst. 19 (2020), 22.

\bibitem{15} L.X. Wang, S.J. Chen, Two solutions for nonhomogeneous
    Klein-Gordon-Maxwell system with sign-changing potential; Elec.
    Diff. Equat. 124 (2018) 1-21.

\bibitem{16} L.P. Xu, H.B. Chen, Existence and multiplicity of
    solutions for nonhomogenous Klein-Gordon-Maxwell equations, Electron. J. Differential Equations 102 (2015) 1-12.

\bibitem{18} Q.Y. Zhang, B. Xu, Multiplicity of solutions for a class
    of    semilinear    Schr\"{o}dinger equations with sign-changing
    potential, J. Math. Anal. Appl. 377 (2011) 834-840.


\bibitem{22} J. Zhang, Solutions to the critical Klein-Gordon-Maxwell
    system with external potential, J. Math. Anal. Appl., 455 (2017) 1152-1177.



\end{thebibliography}
\end{document}